\documentclass[12pt]{article}

\usepackage{amsmath}
\usepackage{amsfonts}

\usepackage{amssymb}

\renewcommand\L{\Lambda}
\newcommand\F{\mathcal F}
\newcommand\G{\mathbf G}
\renewcommand\H{\mathcal H}
\newcommand\J{\mathbf F}

\newcommand\f{\mathbf f}
\newcommand\D{\mathcal D}
\newcommand\LL{\mathcal L}
\renewcommand\l{\lambda}
\newcommand\C{\mathbb{C}}
\newcommand\R{\mathbb{R}}
\newcommand\Z{\mathbb{Z}}
\newcommand\Q{\mathbb{Q}}
\newcommand\QQ{\mathcal Q}
\newcommand\ev{\operatorname{ev}}
\newcommand\ft{\operatorname{ft}}
\newcommand\m{{\mathbf m}}
\newcommand\K{\mathcal K}
\newcommand\I{\mathcal I}
\renewcommand\a{\alpha}
\renewcommand\b{\beta} 
\renewcommand\c{\gamma}
\newcommand\Res{\operatorname{Res}}
\newcommand\w{\wedge}
\newcommand\Lie{\operatorname{Lie}}
\newcommand\p{\partial}
\newcommand\eps{\epsilon}

\title{Gromov--Witten Invariants of Toric Fibrations}

\author{Jeffrey Brown}
\date{}

\begin{document}
\maketitle

\begin{abstract}{We prove a conjecture of Artur Elezi \cite {Elezi} 
in a generalized form suggested by Givental \cite{Givental_MSRI}.
 Namely, our main result relates 
genus-0 Gromov--Witten invariants of a bundle space with such invariants 
of the base, 
provided that the fiber is a toric manifold. When the base is the point, 
a new proof of mirror theorems by A. Givental \cite{Givental_toric}
and H. Iritani \cite{Iritani} for toric manifolds is obtained.}
\end{abstract}    

\section{Formulations}
   
{\bf 1.1. Genus-0 Gromov--Witten invariants.} Given a compact (almost) 
K\"ahler manifold $M$, its {\em genus-0 descendant potential} 
is defined as:
\[ \F_M:=\sum_{n=0}^{\infty}\sum_{D\in MC} \frac{Q^D}{n!}
\int_{[M_{0,n,D}]} \prod_{a=1}^n \sum_{k=0}^{\infty}\ev_a^*(t_k)\psi_a^k.\]
Here $M_{0,n,D}$ stands for the moduli space of degree-$D$ stable maps 
to $M$ of genus-0 holomorphic curves with $n$ marked points, $[M_{0,n,D}]$ ---
its virtual fundamental class, $MC$ --- the Mori cone of $M$, 
i.e. the semigroup of classes in the lattice $H_2(M)$ representable by 
compact holomorphic curves, $Q^D$ --- the element in the {\em Novikov ring} 
(i.e. a power-series completion of the semigroup algebra of the Mori cone) 
representing the degree $D \in MC$ of the stable maps,  $\psi_a$ --- 
the 1st Chern class of the line bundle over $M_{0,n,D}$ formed by cotangent 
lines to the holomorphic curves at the $a$-th marked point, $\ev_a$ --- 
the map $M_{0,n,D}\to M$ defined by the evaluation of stable maps at the 
$a$-th  marked point, $t_k \in H^*(M, \QQ), \ 
k=0,1,2,\dots,$ --- arbitrary cohomology classes of the target manifold $M$
with coefficients in a suitable {\em ground ring} $\QQ$ (for the moment
let it be the rational Novikov ring $\Q [[MC]]$). The explicit inclusion 
of Novikov's variables into the definition of the potential turns out to be
redundant due to the so-called {\em divisor equation} (see 5.1).   
    
Following \cite{Givental_symplectic, Coates-Givental}, one associates to $\F_M$ 
a Lagrangian cone $\LL_M$ in a {\em symplectic loop space} $(\H, \Omega)$. 
Let $H$ denote the cohomology space $H^*(M,\QQ)$, $(\cdot , \cdot)$ 
the Poincare pairing on $H$, and $1\in H$ the unit element. 
Take $\H := H((1/z))$. It consists of Laurent series in one indeterminate 
$1/z$ with vector coefficients. Equip $\H$ with a $\QQ$-valued non-degenerate 
symplectic form 
\[ \Omega (f,g) := \frac{1}{2\pi i} \oint (f(-z),g(z))\ dz .\]
The subspaces $\H_{+}:=H[z]$ and $\H_{-}:=z^{-1}H[[z^{-1}]]$ form a 
Lagrangian polarization of $(\H, \Omega)$, which identifies $\H$ with
$T^*\H_{+}$. To a point $q=q_0+q_1z+q_2z^2+\cdots\ \in \H_{+}$, associate 
a sequence $t=(t_0, t_1, t_2, \dots )$ of elements $t_k \in H^*(X,\QQ)$
according to the {\em dilaton shift} convention:
\[ t_0+t_1z+t_2z^2+\cdots = 1 z + q_0+q_1 z+q_2z^2+\cdots .\]
Define a Lagrangian section $\LL_M$ as the graph of the differential 
of $\F_M$ at the dilaton-shifted point: 
\[ \LL_{M}:=\{ (p,q)\in T^*\H_{+}| p=d_t \F_M \}. \] 
According to general theory of genus-0 Gromov--Witten invariants, the 
section $\LL_M$ considered as a submanifold in $(\H, \Omega)$ is 
(a germ at a dilaton-shifted point of) an {\em overruled Lagrangian cone
with the vertex at the origin}. Here being {\em overruled} means that 
each tangent space $T$ to $\LL_M$ is tangent to $\LL_M$ exactly along
the subspace $zT$ (see \cite{Coates-Givental, Givental_symplectic}).     
This property is invariant under the action of the {\em twisted loop group}
$L^{(2)}GL(H)$. By definition, it consists of those invertible 
Laurent series $W(z)$ with values in $\operatorname{End}(H)$ which preserve
the symplectic form $\Omega$ (i.e. satisfy $W^*(-z)W(z)=1$, where $\,^*$
means ``adjoint'' with respect to the Poincare pairing).   

An overruled Lagrangian cone $\LL_M \subset (\H, \Omega)$ is determined by its
intersection (known as the {J-function}) with the subspace $-1z+z\H_{-}$. 
More precisely, the J-function $\tau \mapsto J(z,\tau)$ is defined as a 
Laurent $1/z$-series with coefficients in $H$ depending on $\tau \in H$ 
and characterized by the property:
\[ J(-z,\tau) = -1z+\tau + O(1/z) \in \LL_M.\]  
Explicitly, for any $\phi\in H$, 
\[ (J(z,\tau),\phi) =(1,\phi)z+(\tau,\phi)+
\sum_{n,D}\frac{Q^D}{n!}\int_{[M_{0,n+1,D}]}
\ev_1^*(\tau)\cdots \ev_n^*(\tau) \frac{\ev_{n+1}^*(\phi)}{z-\psi_{n+1}}.\]

Barannikov \cite{Barannikov}, in a mirror context, constructed a function whose values are obtained as the single intersection points of semi- infinite subspaces in a space of Laurent- series.
      
{\bf 1.2. Toric fibrations.}
K\"ahler toric manifolds can be obtained by symplectic reduction from linear
spaces. 

Let a torus $T^N$ act by diagonal unitary transformations 
on the Hermitian space $\C^N$. Denote by $\mu: \C^N\to \R^N:=\Lie^*(T^N)$ 
the moment map of this action, 
$\mu(z_1,\dots ,z_N)=(|z_1|^2,\dots ,|z_N|^2)$. Let a torus $T^K$ be embedded 
as a subtorus $T^K\subset T^N$. The moment map $\C^N\to \R^K:=\Lie^*(T^K)$
is the composition of $\mu$ with the projection $\R^N\to 
\R^K$ dual to 
the embedding $\Lie T^K\subset \Lie T^N$ of the Lie algebras.
We denote by ${\mathbf m}=(m_{ij}|i=1,\dots , K,\ j=1,\dots, N)$ 
the integer $K\times N$-matrix of this projection. Applying symplectic 
reduction over a chosen value $\omega$ of the moment map, we obtain a 
symplectic toric variety $X=\C^N//_{\omega} T^K$. Since the actions 
of $T^N$ and $T^K$ on $\C^N$ commute, $X$ carries 
a canonical action of $T^N$. 

We will assume that $X$ is non-singular and compact. In fact any compact 
K\"ahler toric manifold $X$ of dimension $N-K$ can be obtained by such 
reduction, with $\R^K$ canonically identified with $H^2(X,\R)$. The 
$T^K$-fibration $({\mathbf m}\circ \mu)^{-1}(\omega) \to X$ endows 
$X$ with $K$ 
tautological $T^N$-equivariant line bundles, whose 1st Chern classes we 
denote by $-p_1,\dots, -p_K$. They represent a basis in $\R^K$ of integer lattice vectors, 
and generate the algebra $H^*(X)$.    

Let $B$ be any K\"ahler manifold, $L_1,\dots ,
L_N$ line bundles over $B$, and $\Lambda_{j}= c_{1}(L_{j}^{*}),\ j=1, \cdots, N$.
In the vector bundle $\oplus L_j$ with the structure group $T^N$, replace the
fiber with the toric $T^N$-space $X$. We obtain a toric fibration $\pi: 
E\to B$ of K\"ahler manifolds. It carries a canonical fiberwise action 
of $T^N$. The total space $E$ is endowed with $K$ tautological line bundles
whose 1st Chern classes we denote by $-P_1,\dots ,-P_K$.
They restrict to the fibers to $-p_1,\dots, -p_K$, and generate $H^*(E)$ 
as an algebra over $H^*(B)$.   
 
To a degree $\D \in H_2(E)$ of holomorphic curves in $E$, we associate 
the degree $D:= \pi_*(\D) \in H_2(B)$ of its projection to the base, and
the degrees $d_i:=P_i(\D),\ i=1,\dots , K$, with respect to the classes $P_i$. 
In the Novikov ring of $E$, we will represent $\D$ by the monomial 
$q^d Q^D$, where $q^d=q_1^{d_1}\cdots q_K^{d_K}$, and $Q^D$ represents $D$ in 
the Novikov ring of $B$. 

In the formulation below we use the following notation: 

\begin{align*} t=(t_1,\dots, t_K),\ \ Pt=\sum_{i=1}^K P_it_i,\ \ 
dt=\sum_{i=1}^K d_it_i, \\
\ U_j=\sum_{i=1}^K P_im_{ij}-\L_j, \ \ 
U_j(\D)=\sum_{i=1}^K d_im_{ij}-\L_j(D).\end{align*}

{\bf 1.3. Main results.} 

{\tt Theorem 1.} {\em Decompose the J-function of the overruled 
Lagrangian cone $\LL_B$, corresponding to the base $B$ of a toric
fibration $E\to B$, according to the degrees of curves:
\[ J(z,\tau) =\sum_{D\in MC(B)} J_D(z,\tau) Q^D, \]
and introduce the {\bf \em hypergeometric modification} 
\[ I_E(z,t,\tau,q,Q):=e^{Pt/z}
\sum_{d\in \Z^K, D\in MC(B)} \frac{J_D(z,\tau) Q^D q^de^{dt}}
{\prod_{j=1}^N \prod_{m=1}^{U_j(\D)} (U_j+mz)}.\]
Then for all $(t,\tau)$, the series $I_E(-z)$ lies in the overruled 
Lagrangian cone $\LL_{E}$ corresponding to the total space $E$.}

The products in the denominator are interpreted as ratios of the values of the
Gamma-function:
\[ \prod_{m=1}^n (U+mz) := \prod_{m=-\infty}^n (U+mz)\ \ / \prod_{m=-\infty}^0
(U+mz).\]
Thus, when $n<0$, we obtain a product in the numerator instead.

In the special case when $E = \operatorname{proj} \left(\oplus_j L_j\right)$ 
is a projective fibration over $B$, the algebra $H^*(E)$ is generated 
over $H^*(B)$ by one generator $P$ satisfying the relation 
$(P-\L_1)\cdots (P-\L_N)=0$.    

{\tt Corollary 1.} {\em The overruled Lagrangian cone $\LL_E$ of the 
projective fibration contains $I(-z)$, where 
\[ I(z,t,\tau,q,Q):=e^{Pt/z}\! \sum_{D\in MC(B)}\! J_D(z,\tau) Q^D 
\sum_{d\in \Z} \frac{q^d e^{dt}}{\prod_{j=1}^N \prod_{m=1}^{d-\L_j(D)}
(P-\L_j+mz)}.\]}  

Note that the summation range $d\in \Z$ actually reduces to $d\geq 
\min_j \L_j(D)$ since otherwise the numerator contains  
$\prod_j (P-\L_j)=0$. We will see later that a similar phenomenon takes place
in the situation of general toric fibrations. As a result, the effective
summation range in the series $I_E$ stays within the Mori cone of $E$.

{\tt Corollary 2} (Elezi's conjecture \cite{Elezi}).{\em 
Taking $\L_1=0$, and assuming that all $\ L_j^{*}$ and $c_1(T_B)- \sum_j \L_j$      
are {\em nef}, we have: $I = z+\tau+Pt +O(z^{-1})$, i.e.
the series $I$ represents the J-function of the projective fibration 
$E\to B$ at the points $\tau+Pt \in H^*(E, \QQ)$.}

When $B=pt$, we have $J_B=ze^{\tau/z}$, leading to Iritani's mirror 
theorem for arbitrary toric manifolds. 
 
{\tt Corollary 3} (Iritani's theorem \cite{Iritani}.)
{\em For all values of $(\tau, t)$, the series $I_X(-z)$, where  
\[ I_X(z,t,\tau,q):=z \, e^{(\tau+Pt)/z} \sum_d \frac{q^de^{dt}}
{\prod_{j=1}^N \prod_{m=1}^{\sum_i d_im_{ij}} (\sum_i P_im_{ij}+mz)}, \]
lies on the overruled Lagrangian cone $\LL_X$ of the toric manifold $X$.}

{\tt Corollary 4} (Givental's theorem \cite{Givental_toric}).
{\em When the toric manifold $X$ is Fano, then the series $I_X$ 
represents the J-function of $X$ at the points $\tau \oplus Pt \in 
H^0(X)\oplus H^2(X)$.} 
     
{\bf 1.4. Remarks.}

Although the results are stated above in the setting of K\"ahler manifolds, they extend without complications to the general
setting of symplectic toric fibrations and almost K\"ahler structures. 

Furthermore, in all formulations one may assume that cohomology groups 
are {\em equivariant} with respect to the fiberwise action of the torus 
$T^N$ on the toric fibration $E\to B$. Respectively, all Gromov--Witten 
invariants become equivariant, taking values in the coefficient ring 
$H^*(BT^N, \Q)=\Q [\l_1,\dots , \l_N]$ of the equivariant cohomology 
theory. In this case, the classes $P_i$ and $\L_j$ are understood as 
$T^N$-equivariant 1st Chern classes of the respective line bundles, and are
invertible in the field of fractions of the coefficient ring. 
In fact Theorem 1 follows in the limit $\l=0$ from its equivariant 
counterpart. 

To prove the equivariant version of Theorem 1, we will first 
show in Section 2 that the equivariant counterpart of the 
overruled Lagrangian cone $\LL_E$ is the solution set of a certain 
recursion relation. This is an unpublished result of A. Givental, 
and an easy special case of the general fixed point localization 
formula for Gromov--Witten invariants in the case of non-isolated fixed 
points \cite{Brown et al}. 

Cohomology classes entering in the definition of $I_E$ have equivariant counterparts, so we may consider $I_E$ as taking values in equivariant cohomology.  The equivariant version of the hypergeometric modification series $I_E$ 
has essential singularity at $z=0$ and simple poles at $z\neq 0$.  To prove that the series satisfies the recursion relation, one needs to show that:
(i) applying a certain linear transformation removes the essential 
singularity at $z=0$, and (ii) residues at the simple poles are controlled 
recursively. In Section 3, we show (ii) by decomposing terms of the series
$I_E$ into elementary fractions in a straightforward way. The task (i) relies 
on properties of oscillating integrals arising in mirror theory. It is 
accomplished in Section 5 in a way resembling the proof \cite{Coates-Givental}
of Quantum Lefschetz Theorem. This is preceded by a general discussion 
in Section 4 of asymptotics of oscillating integrals.       

\section{Localization}

{\bf 2.1. Fixed sections.} Let $X=\C^N//_{\omega}T^K$ be a compact toric 
manifold as in 1.2. For $X$ to be non-singular, it is necessary that 
$\omega \in \R^K$ is a regular value of the moment map $\m \circ \mu$. 
The image of the moment map is a {\em picture} of the 1st orthant $\R_{+}^N$ 
``drawn'' (by means of the projection $\m$) in $\R^K$. The fiber 
$\m^{-1}(\omega)\subset \R_{+}^N$ is the momentum polyhedron of the torus 
$T^N/T^K$ action on $X$. The vertices of the momentum
polyhedron represent fixed points of the torus action on $X$.
Each vertex corresponds to a $K$-dimensional face of the 1st orthant whose
picture contains $\omega$. We will label the fixed point by the multi-index 
$\a =\{ j_1<\dots < j_K \}$ specifying the coordinates of 
the corresponding $K$-dimensional face.  
\footnote{Pictures in $\R^K=H^2(X,\R)$ of such $K$-dimensional faces 
contain, together with $\omega$, its connected component $\K$ in the regular 
value locus of the moment map. $\K$ coincides with the K\"ahler cone of $X$.} 

Consider now a toric fibration $\pi: E\to B$ with the fiber $X$ (as in 1.2).
Fixed points of the torus $T:=T^N$ acting fiberwise on $E$ form sections 
$\a: B\to E$, one for each fixed point $\a \in X^T$ of the torus action on 
the fiber. Torus-equivariant intersection theory on the total space $E$ 
of the fibration can be completely characterized in terms of intersection 
theory on the base $B$ by the following elegant residue formula describing 
(via fixed point localization) the push-forward to $H^*_T(B)=H^*(B, H^*(BT^N))$
 of a $T^N$-equivariant cohomology class $f \in H^*_T(B)[P_1,\dots , P_K]$:
\[ \pi_* f = \sum_{\a \in X^T} \Res_{\a} f(P) 
\frac{dP_1\w \cdots \w dP_K} {U_1(P)\cdots U_N(P)}.\]
The factors $U_j=\sum_i P_i m_{ij}-\L_j$, $j=1,\dots, N$, can be 
interpreted as Poincare-duals of the torus-invariant divisors represented 
by hyperplane faces of the momentum polyhedron $\m^{-1}(\omega)$. 
Here $\Res_{\a}$ refers to the residue of the $K$-form at the pole
$U_{j_1}=\cdots =U_{j_K}=0$ corresponding to the fixed point 
$\a=\{ j_1,\dots , j_K\}$, i.e. at the point $P=P^{\a}$ determined by 
\[ \sum_i P^{\a}_i m_{ij}=\L_j, \ \forall j\in \a.\]
The formula for the push-forward uses the wedge product symbol in a non- standard way.  Namely, write 
\[dP_1 \wedge \dots \wedge dP_k= { {dP_1 \wedge \dots \wedge dP_k}\over {dU_{j_{1}} \wedge \dots \wedge dU_{j_{K}}}}dU_{j_{1}} \wedge \dots \wedge dU_{j_{K}}\]
where the ratio of $K$- forms is equal to $det^{-1} (m_{i, j_{s}})$, which is $\pm 1$ for smooth toric fibers.  The wedge symbol, as it is used in the formula for the push-forward, is an instruction to compute residue integrals with the $dP_i$'s reordered according to the exterior algebra so as to offset this sign.

We note that: (i) the normal bundle to the fixed point section $\a$ is the 
sum of $N-K$ line bundles with the 1st Chern classes 
\[ \a^*U_j =\sum_iP^{\a}_im_{ij}-\L_j, \ \text{where}\ j\notin \a,\]
and (ii) a point $D \in MC(B)$, lifted to $E$ by this section, is represented
in the Novikov ring of $E$ by the monomial 
$Q^Dq^{P^{\a}(D)}=Q^D q_1^{P_1^{\a}(D)}\cdots 
q_K^{P_K^{\a}(D)}$.

{\bf 2.2. The cone $\LL_E$.} The overruled Lagrangian cone $\LL_E$ in the
torus-equivariant genus-0 Gromov--Witten theory of the total 
space $E$ of the toric fibration lies in the appropriate symplectic 
loop space $(\H, \Omega)$. 
The space is actually a module over the ground ring $\QQ$, which we 
currently take to be the Novikov ring of $E$ tensored with the field 
of fractions $\Q (\l )$ of $H^*(BT^N)$. Pending further completions, $\H$ 
consists of Laurent series in $1/z$ with coefficients in $H=H^*(E, \QQ)$.  
A point in the cone can be written as   
\[ \J (-z, t) = -1z+t(z)+\sum_{n,D,d}\frac{Q^Dq^d}{n!}(\ev_1)_*
\left[ \frac{1}{-z-\psi_1} 
\prod_{i=2}^{n+1} (\ev_i^*t)(\psi_i)\right],\]
where $(\ev_1)_*$ denotes the virtual push-forward by the evaluation
map $\ev_1: E_{0,n+1,\D}\to E$, and $t(z)=\sum_{k=0}^{\infty}t_kz^k$
is a polynomial with arbitrary coefficients $t_k\in H$. 

Denote by $\J^{\a}:=\a^*\J$ restrictions of $\J$ (considered as a cohomology 
class of $E$) to the fixed point sections $\a$. The series $\J^{\a}$ lie 
in the space of Laurent series in $1/z$ with coefficients in $H^*(B, \QQ)$.
In terms of the push-forwards $\a_*: H^*(B,\QQ) \to H^*(E,\QQ)$ by the 
sections and their normal Euler classes $e^{\a}$, we have:
\[ \J = \sum_{\a\in X^T} \a_*\left(\frac{\J^{\a}}{e^{\a}}\right), \ \ 
\text{where}\ 
e^{\a}=\prod_{j\notin \a} U_j(P^{\a}).\] 

{\bf 2.3. Twisted Gromov--Witten invariants.} Consider the base $B$ of the
toric fibration embedded in the total space $E$ as a fixed section 
$\a: B\to E$. Torus-equivariant Gromov--Witten invariants of a {\em 
neighborhood} of this section can be defined via fixed point localization 
as certain intersection indices in moduli spaces of stable maps to the 
fixed locus $\a(B)$. They coincide with such invariants
of $B$ {\em twisted} (in the sense of \cite{Coates-Givental}) by the normal
bundle $N^{\a}$ of the fixed section in $E$. More specifically, the genus-$0$
descendant potential of the twisted theory is defined by the formula:\footnote{
Here $Q^Dq^{P^{\a}(D)}$ represents degree-$D$ curves of $B$ considered
as curves in $\a(B)\subset E$.}   
\[ \F_{B,N^{\a}}:=\sum_{n,D} \frac{Q^Dq^{P^{\a}(D)}}{n!}\int_{[B_{0,n,D}]}
Euler_T^{-1}(N^{\a}_{0,n,D})
\prod_{a=1}^n \sum_{k=0}^{\infty}\ev_a^*(t_k)\psi_a^k.\]      
Here $Euler^{-1}_T$ is the inverse $T$-equivariant Euler class of complex
vector bundles\footnote{In \cite{Coates-Givental}, twisting by arbitrary
invertible multiplicative characteristic classes is allowed.},
and $N^{\a}_{0,n,D}:=(\ft_{n+1})_*\ev_{n+1}^*N^{\a}$ denotes 
the virtual vector bundle over $B_{0,n,D}$ obtained as the K-theoretic 
push-forward along the family of curves 
$\ft_{n+1}: B_{0,n+1,D} \to B_{0,n,D}$ of the bundle $N^{\a}$
pulled-back from $B$ by $\ev_{n+1}: B_{0,n+1,D}\to B$.

Let $\LL^{\a}$ be the overruled Lagrangian cone corresponding 
to the twisted theory. The cone lies in the symplectic loop space 
$(\H^{\a},\Omega^{\a})$ constructed using the twisted Poincare pairing
$(a,b)^{\a}=\int_B Euler^{-1}_T(N^{\a}) ab$ on $H^{\a}:=H^*(B,\QQ)$.
Let $T$ denote a tangent space to the
cone $\LL^{\a}$ at a point $\f$. The same space $T$ is tangent to
$\LL^{\a}$ everywhere along $zT$. Let $-z+u^{\a}+O(1/z)$ be the point
of the J-function of this cone that lies in $zT$. Here $u^{\a}$ depends on 
$\f\in \LL^{\a}$, and is an element of $H^{\a}$.   

To each vector $w\in H^{\a}$, associate the vector in $T$ that projects to
$w$ along $\H_{-}^{\a}$. The operator thus defined is an element of 
the loop group $LGL(H^{\a})$, and is represented by an operator Laurent series
of the form $1+O(1/z)$. The inverse element, which we denote 
by $S_{u^\a}(-z)$, is characterized as the operator $1/z$-series which 
transforms $T$ to $\H^{\a}_{+}$. 
From the fact that $T$ is Lagrangian, it follows that 
$S^*_{u^\a}(z)S_{u^\a}(-z)=1$ (i.e. $S_{u^{a}}$ lies in the twisted 
loop group). We conclude that for every point $\f\in \LL^{\a}$ there exists a 
unique $u^{\a}\in H^{\a}$ such that $S_{u^{\a}} \f  \in  z\H^{\a}_{+}$.

{\bf 2.4. Fixed stable maps.} The general description of genus $0$ 
stable maps $\Sigma \to M$ whose equivalence class is fixed by the action 
of a torus $T$ on the target space goes back to Kontsevich's work 
\cite{Kontsevich}. According to it, each irreducible component of 
the curve $\Sigma$ must be mapped onto an orbit of dimension $0$ or $1$ of the 
complexified torus $\,^{\C}T$. A $1$-dimensional orbit is a projective 
line connecting two $0$-dimensional orbits (i.e. fixed points). An irreducible 
component of $\Sigma$ mapped onto such orbit with degree $k$ must have 
ramifications of degree $k$ over the fixed points. We will call such 
irreducible components {\em legs} (of {\em multiplicity} $k$). 
Removing all legs from 
$\Sigma$ leaves a forest of rational curves mapped to the fixed point locus 
$M^T$ in the target space. Integration over fixed point components in the 
moduli spaces of stable maps reduces therefore to evaluation of certain 
twisted Gromov--Witten invariants of $M^T$. 

In our situation, $1$-dimensional orbits of the torus $\,^{\C}T$ 
in the toric manifold $X$ correspond to edges of the momentum polyhedron. 
If two vertices of the polyhedron are connected by an edge, then there 
is exactly one $1$-dimensional orbit connecting two corresponding fixed points
(say, $\a$ and $\b$) in $X$. Respectively, each fiber of the toric 
fibration $E\to B$ contains a copy of this $1$-dimensional orbit, connecting
the fixed point sections $\a$ and $\b$.   

Two fixed points $\a$ and $\b$ are connected by a $1$-dimensional orbit exactly
when the union $\a\cup \b$ of the multi-indices $\a$ and $\b$ has cardinality 
$k+1$. The orbit itself is a toric $\C P^1=\C^{k+1}//_{\omega}T^k$ obtained 
by symplectic reduction from the face of the 1st orthant whose coordinates
have indices from $\a\cup\b$. From this, one can easily derive the following
relations (see \cite{Givental_toric}). Denote by $j_{\pm}(\a,\b)$ the indices 
such that 
$j_{+}(\a,\b) \in \b - \a$, and $j_{-}(\a,\b)\in \a-\b$, by $\chi_{\a,\b}$
the equivariant $1$-st Chern class of the line bundle over $B$ formed 
by the tangent lines to the $1$-dimensional orbit at the fixed points $\a$, and
by $d_{\a,\b}$ the degree of the $1$-dimensional orbit as a rational curve 
in the fiber $X$. Then for a given fixed point $\a$, any index 
$j_{+}\notin \a$ can play the role
of $j_{+}(\a,\b)$, while $\b$ and $j_{-}(\a,\b)$ are uniquely determined
by it. By fixed point localization on $\C P^1$, we find:
\[ d_{\a,\b}=\frac{P^{\a}-P^{\b}}{\chi_{\a,\b}}, \ \
U_j(d_{\a,\b})=\frac{\a^*U_j-\b^*U_j}{\chi_{\a,\b}},\]
where $\chi_{\a,\b}=\a^*U_{j_{+}(\a,\b)}=-\b^*U_{j_{-}(\a,\b)}$
It follows that
\[ U_{j_{\pm}(\a,\b)}(d_{\a,\b})=1 \ \text{and}\ 
\forall\ j\in \a\cap\b, \ \a^*U_j=\b^*U_j=0, \ U_j(d_{\a,\b})=0.\]

{\bf 2.5. Recursion.} 
In general, the value of a $T$-equivariant cohomology class
$f$ on the invariant fundamental class of a manifold (or orbifold) $M$ is 
computed as 
\[ \int_{M} f = \int_{M^T} \frac{i^*f}{Euler_T({\mathcal N})},\]
where $i: M^T\to M$ is the embedding of the fixed point locus, 
${\mathcal N}$ is the normal bundle to $M^T$ in $M$, and $Euler_T$ is the
$T$-equivariant Euler class. The use of fixed point localization 
in application to 
integrals over {\em virtual} fundamental classes of moduli spaces 
of stable maps has been justified 
by Graber--Pandharipande \cite{Graber-Pandharipande}.  
Our nearest goal is to characterize points of the cone $\LL_E$ by 
a recursion relation which comes from fixed point localization in moduli 
spaces $E_{0,n+2,\D}$.  

Let $\J\in \LL_E$ and $\J^{\a}:=\a^*\J$, i.e.  
\[ \J^{\a} (-z, t) = -1z+\a^*t(z)+\a^*\sum_{n,D,d}\frac{Q^Dq^d}{n!}(\ev_1)_*
\left[ \frac{1}{-z-\psi_1} 
\prod_{i=2}^{n+1} (\ev_i^*t)(\psi_i)\right].\]
We evaluate the sum via fixed point localization, and notice first of all, 
that a torus-fixed stable map $\Sigma \to E$ does not contribute to 
$\J^{\a}$ unless the 1st marked point lands in the fixed section $\a$. 
When it does, there are two possibilities: the marked point can belong 
to a {\em leg} (see 2.4), or to a tree $C \subset \Sigma$ 
of rational components mapped to the locus 
$\a(B) \subset E$ of the fixed point section $\a: B\to E$.  

Examine the first possibility. The leg carrying the 1st marked point 
is a ramified cover of multiplicity $k>0$ of a 1-dimensional orbit of the 
torus $\ ^{\C}T$ which lies in a fiber $X$ of the toric fibration $E\to B$,
connects the fixed point $\a$ with another fixed point $\b$, and has the 
degree $d_{\a,\b}$ considered as a curve in $X$.
Contributions of all stable maps of this type to $\J^{\a}$ via fixed 
point localization can be represented in the form:\footnote{Here and later
$\b=\b(j_{+})$ where $j_{+}$ runs all indices not in $\a$, 
as explained in 2.4.}
\[ \sum_{\b}\sum_{k>0} \frac{q^{kd_{\a,\b}}\ 
Euler^{-1}_T({\mathcal N}_{\a,\b}(k))}{k(-z+\chi_{\a,\b}/k)} 
\ \J^{\b}(-\frac{\chi_{\a,\b}}{k}).\]  
Here the factor $-z+\chi_{\a,\b}/k$ is the specialization of $-z-\psi_1$
to the fixed point component (namely, $-\chi_{\a,\b}/k$ is the equivariant
1st Chern class of the line bundle over $B$ formed by the cotangent lines 
to the leg at the 1st marked point). The symbol ${\mathcal N}_{\a,\b}(k)$ 
denotes the virtual bundle over $B$ whose fibers describe deformation modes
of the 1st leg (in the direction normal to the fixed point locus in the moduli
spaces of stable maps to $E$). The occurrence of the Euler class of this 
bundle in the denominator is due to the general structure of localization 
formulas. One can easily compute this Euler class explicitly:
\[  Euler_T({\mathcal N}_{\a,\b}(k))= \prod_{m=1}^{k-1} 
(\a^*U_{j_{+}(\a,\b)}-m\frac{\chi_{\a,\b}}{k})
\prod_{j\notin \b} \prod_{m=1}^{kU_j(d_{\a,\b})} 
(\a^*U_j-m\frac{\chi_{\a,\b}}{k}).\]
The normal bundle to the fixed point locus in the moduli space 
contains the smoothing mode of the curve $\Sigma$ at the node where 
the leg and the rest of the curve, $\Sigma'$, connect. This is a line bundle 
with the 1st Chern class $-\psi-\chi_{\a,\b}/k$, where $\psi$ is the 1st
Chern class of the universal cotangent line bundle to the curves $\Sigma'$
at the node. Considering the node as the 1st marked point of the curve 
$\Sigma'$, we can therefore represent the sum of all contributions of moduli
spaces of curves $\Sigma'$ as $\J^{\b}(-\chi_{\a,\b}/k)$. The extra factor
$k$ in the denominator is due to the cyclic symmetry of order $k$ of the leg,
which affects the orbifold's fundamental class this way.        

Now put
\[ t^{\a}(z):=\a^*t(z)+\sum_{\b}\sum_{k>0} \frac{q^{kd_{\a,\b}}\ 
Euler^{-1}_T({\mathcal N}_{\a,\b}(k))}{k(-z+\chi_{\a,\b}/k)} 
\ \J^{\b}(-\frac{\chi_{\a,\b}}{k}),\]
and examine the second possibility, when the 1st marked point of the
curve $\Sigma$ lies on a tree $C$ mapped to $\a(B)$. Fixed point components
of the moduli spaces of stable maps to $E$ are products of moduli spaces of
stable maps to the fixed point sections, and such moduli of the maps 
$C\to \a(B)$ is one of the factors. Integrating over this factor last, we
represent the contribution of each fixed point component as a genus-0
Gromov--Witten invariant of $B$ twisted by the normal bundle $N^{\a}$.
Among the marked points of the curves $C$, one is the 1st marked point of
$\Sigma$; it carries the input $1/(-z-\psi_1)$. Every other marked point
could be either a marked point of $\Sigma$ which happens to lie in $C$, 
or a node where a connected component of $\Sigma - C$ is attached to $C$ 
by a leg. The input at the marked point 
of this is obtained by adding to $\a^*t(\psi_i), (i>1)$ the sum of 
fixed point localization contributions over all possibilities for the 
connected component of $\Sigma-C$. The total input coincides with 
what is denoted above by $t^{\a}(\psi_i)$. Thus we have:
\begin{align*} \J^{\a}(-z)&=-z+t^{\a}(z)+ \\ &\sum_{n,D}
\frac{Q^Dq^{P^{\a}(D)}}{n!}\ (\ev_1)_*\left[ \frac{
Euler_T^{-1}({\mathcal N}_{0,n+1,D})}{-z-\psi_1}
\prod_{i=2}^{n+1}t^{\a}(\psi_i)\right],\end{align*}
where $(\ev_1)_*$ is the push-forward in the $N^{a}$-twisted 
Gromov--Witten theory $B$ by the map $\ev_1: B_{0,n+1,D}\to B$. 
We conclude that {\em $\J^{\a}(-z)$ lies in the overruled Lagrangian cone
$\LL^{\a}$ of the twisted Gromov--Witten theory of $\a(B) \subset E$.} 

Let us examine analytical properties of the expression for $\J^{\a}$ as a 
function of $z$. Since the class $\psi_1$ in the sum is nilpotent, each summand
with a fixed $D$ and $n$ is polynomial in $1/z$. When $D$ is fixed, but $n$ 
grows, the degree of the polynomial can grow too. In fact, employing 
dimensional arguments and the string equation, one can see that for a fixed
$D$ the sum over $n$ is a finite linear combination of functions of the form
$z^{-k}e^{c/z}$ with positive $k$ and non-zero constant $c$. 
Thus the whole sum is a $Q$-series whose coefficients are meromorphic 
functions with {\em essential 
singularities} at $z=0$ and no other singularities (including $z=\infty$).    
The term $t^{\a}(z)$, in the contrary, is a $q$-series, whose coefficients
have simple poles at $z=\chi_{\a,\b}/k$, and (from the summand $\a^*t(z)$) 
a pole at $z=\infty$ of any order $\geq 0$. Thus, {\em $\J^{\a}(z)$ are
power series in the Novikov variables $Q,q$ which have coefficients 
meromorphic in $z$, with an essential singularity at $z=0$, 
finite order pole at $z=\infty$ and simple poles at $z=-\chi_{\a,\b}/k$, such
that the residues at the simple poles satisfy the recursion relation:
\[ \operatorname{Res}_{z=-\frac{\chi_{\a,\b}}{k}} \J^{\a}(z)\ dkz = 
\frac{q^{kd_{\a,\b}}}{Euler_T({\mathcal N}_{\a,\b}(k))}\ 
\J^{\b}(-\frac{\chi_{\a,\b}}{k}).\]} 
\indent Note that each elementary fraction $(z+\chi_{\a,\b}/k)^{-1}$ 
can be expanded into a Laurent series in two ways: inside or outside the circle
$|z|=|\chi_{\a,\b}|/k$. Expanding all the elementary fractions as $1/z$-series 
renders $\J^{\a}(-z)$ as a Laurent series in $1/z$ in a way it occurs as a
component of $\J (-z) \in \LL_E$. Expanding all the elementary fractions as
$z$-series renders $\J^{\a}(-z)$ as a Laurent series in $1/z$ in another 
way, namely the way it occurs as a point on the cone $\LL^{\a}$. 
Note that the polynomial truncation $t^{\a}=[\J^{\a}(z)]_{+}$ of the latter 
Laurent series lies not in $H^{\a}[z]$ per se, but in a certain completion
of it. Namely, it is a power series in $1/\chi_{\a,\b}$ with coefficients
polynomial in $z$, or equivalently becomes a polynomial in $z$ when reduced 
modulo any power of $1/\chi_{\a,\b}$. We will subsequently assume that 
{\em the ground ring is suitably localized to include inverse powers of 
$\chi_{\a,\b}$,} and that {\em $\LL^{\a} \subset \H^{\a}$ refers to the cone 
of the twisted Gromov--Witten theory thus completed.} With these 
interpretations in mind, we state the result of this section.

{\tt Theorem 2.} {\em Points $\{ \J^{\a}(-z) \}$ of the overruled  Lagrangian 
cone $\LL_E$ are characterized by the following conditions: 

(i) $\J^{\a}(-z) \in \LL^{\a}$, 

(ii) $\J^{\a}$ are power series in the Novikov variables
$Q,q$ whose coefficients are analytic functions of $z$ with  
essential singularities at $z=0$, finite order poles at $z=\infty$, 
simple poles at $z=\chi_{\a,\b}/k,\ k=1,2,3,\dots $, 
and such that the residues 
at the simple poles satisfy the recursion relations:
\[ \operatorname{Res}_{z=-\frac{\chi_{\a,\b}}{k}} \J^{\a}(z)\ dkz =    
\frac{q^{kd_{\a,\b}}}{Euler_T({\mathcal N}_{\a,\b}(k))}\ 
\J^{\b}(-\frac{\chi_{\a,\b}}{k}),\]
where  
\[ Euler_T({\mathcal N}_{\a,\b}(k))= \prod_{m=1}^{k-1} 
(\a^*U_{j_{+}(\a,\b)}-m\frac{\chi_{\a,\b}}{k})
\prod_{j\notin \b} \prod_{m=1}^{kU_j(d_{\a,\b})} 
(\a^*U_j-m\frac{\chi_{\a,\b}}{k}).\]}

We have established that every point on $\LL_E$ satisfies (i) and (ii).
Let us prove now that if series $\{ \J^{\a}\}$ satisfy these conditions, 
then they represent a point in $\LL_E$. Indeed, since $\J^{\a}(-z) 
\in \LL^{\a}$, there exists a unique $u^{\a}\in H^{\a}$ such that 
$\G^{\a}(z):=S_{u^{\a}}(-z) \J^{\a}(z) \in z\H^{\a}_{+}$ (see 2.3).    
Combining this with the property (ii), we conclude that as a function of $z$,
$\G^{\a}$ satisfies the following recursion relation:
\begin{align*} (\ast)\ \ \G^{\a}(z)&=q^{\a}(z)+ \\ 
&\sum_{\b, k} S_{u^{\a}}(\frac{\chi_{\a,\b}}{k}) 
\frac{q^{k d_{\a,\b}} Euler_T^{-1}({\mathcal N}_{\a,\b}(k))}{kz+\chi_{\a,\b}} 
S^*_{u^{\b}} (-\frac{\chi_{\a,\b}}{k}) \G^{\b}(-\frac{\chi_{\a,\b}}{k}),
\end{align*}
where $q^{\a}$ are $(Q,q)$-series with coefficients polynomial in $z$.
Indeed, for $m>0$,
\[ \frac{z^{-m}-(-\chi_{\a,\b}/k)^{-m}}{kz+\chi_{\a,\b}} = O(\frac{1}{z}).\]
Since $S_{u^{\a}}(-z)=1+O(1/z)$ is a $1/z$-series, we see that the difference
between $\G^{\a}(z)$ and the R.H.S. of the recursion relation is $O(1/z)$, and 
hence vanishes, because both sides lie in $\H^{\a}_{+}$. 

Given arbitrary $\{ q^{\a} \}$, and $\{ u^{\a}\}$, a solution 
$\{ \G^{\a} \}$ to the system of recursion relations $(\ast)$ is computed 
by successive $q$-adic approximations, and is therefore unique. 
The set of corresponding $\J^{\a} \in \LL^{\a}$ is reconstructed by the 
application of $S_{u^{\a}}^{-1}(-z)$. In particular, $\a^*t(z)$ are related to
$q^{\a}(z)$ by  
\[ z+\a^*t(-z) = \left[ S_{u^{\a}}^{-1}(-z) q^{\a}(z) \right]_{+}.\]
It remains to show that the values $\{ u^{\a} \}$ are unambiguously determined 
by $\{ \a^*t \}$ in view of the additional constraint that 
$\G^{\a} \in z\H^{\a}_{+}$ (rather than $\H^{\a}_{+}$). Indeed, if such 
uniqueness is established, we conclude that $\{ \J^{\a}\}$ coincide with the 
components of the point $\J$ on the cone $\LL_{E}$ which corresponds to the
Gromov--Witten invariants of $E$ with the inputs $t$ at the marked 
points. 

To verify the required uniqueness, consider first the same problem 
{\em classically}, i.e. modulo Novikov variables $Q$ and $q$. Then recursion 
relation $(\ast)$ degenerates into $\G^{\a}=q^{\a}$, the S-matrices turn into
$S_{u^{\a}}(-z)=e^{-u^{\a}/z}$ (where $u^{\a}\in H$), so that we have:
\[ q^{\a}(-z)=\left[ e^{u^{\a}/z} (-z+\a^*t(z))\right]_{+}. \]
The additional constraints $q^{\a}(0)=0$ assume the same form of the universal
fixed point equation
\[ u =\sum_{m=0}^{\infty} t_m \frac{u^m}{m!},\]
where $t_m$ are coefficients of $\a^*t =\sum t_m z^m$. 
The fixed point equation has a unique formal solution $u=u(t_0,t_1,t_2,\dots)$
on the space of polynomials. Using the formal Inverse Function Theorem, we 
conclude that the values of $u^{\a}$ can be uniquely found by successive 
$(Q,q)$-adic approximations from the relations between $\{ \a^*t \}$ and 
$\{ q^{\a} \}$, the recursion relations $(\ast)$, and the additional 
constraints $\G^{\a}\in z\H^{\a}_{+}$.  
  
\section{Recursion}

To prove the equivariant version of Theorem 1, it suffices to show that 
$\J=I_{E}$ satisfies conditions (i) and (ii) of Theorem 2. The hypergeometric
modification $I_E$ is a $(q,Q)$-series whose coefficients have simple poles
at $z=-\a^*U_j/k$, finite order poles at $z=\infty$, and essential 
singularities at $z=0$. Thus we need to show
that: (i) $\a^*I_E \in \LL^{\a}$, and (ii) residues at the simple poles 
satisfy the recursion relation of Theorem 2. We postpone (i) until Section 5, 
and deal with (ii) here by computing the residues explicitly. 

We have:
\begin{align*} \J^{\a}(z):=\a^* I_E(z) =
e^{P^{\a}t/z} 
\sum_{D}\sum_{d'\in {\mathbb Z}^K} \frac{J_D(z,\tau) Q^D q^{d'}e^{d't}}
{\prod_{j=1}^N \prod_{m=1}^{\sum_i d'_i m_{ij}-\L_j(D)} (\a^*U_j+mz)}.
\end{align*} 
It will be convenient to put $d_i:=d'_i-P_i^{\a}(D)$,       and use that 
\[  \sum_i P^{\a}_i(D)m_{ij}-\L_j(D) = \a^*U_j(D) (=0 \ \forall j\in \a),
\] 
to obtain:
\begin{align*} \J^{\a}(z)=e^{P^{\a}t/z} \sum_D\sum_{d \in \Z^K} 
\frac{J_D(z,\tau) (Q^Dq^{P^{\a}(D)}) q^d e^{dt} e^{P^{\a}(D)t}}
{\prod_{j\in \a} \prod_{m=1}^{U_j(d)} (mz)\ \prod_{j\notin \a} 
\prod_{m=1}^{U_j(d)+\a^*U_j(D)} (\a^*U_j+mz)}.
\end{align*}  
We see that if $U_j(d) < 0$ for some $j\in\a$, then the term contains
a factor $(0z)$ in the numerator. Thus, the effective summation range is
over those $d$ for which $U_j(d) \geq 0$ for all $j\in \a$. 
For each fixed point $\a$, this range lies in the Mori cone of the fiber
$X$ of our toric fibration (because the K\"ahler cone of $X$ lies in 
the simplicial cone spanned by $\{ U_j | j\in \a\}$).  
Since the monomials $Q^Dq^{P^{\a}(D)}$ represent
degrees of holomorphic curves in the fixed section $\a: B\to E$, we conclude
that all series $\a^*I_E$ are supported in the Mori cone of $E$. 
The same remains true for the non-equivariant limit of $I_E$ (as we promised 
in 1.3).

Non-zero poles of $\J^{\a}$ correspond to
the choice of a factor $\a^*U_j+kz$ with $j\notin \a$ and $k>0$.
Given a choice, we put $\a \cup \b = \a \cup \{ j \}$. As we mentioned in 2.4,
this determines $\b$, $d_{\a,\b}$, and $j_{\pm}(\a,\b)$ such that 
$j_{+}(\a,\b)=j$. To single out the contribution of the elementary fraction 
$(\a^*U_{j_{+}(\a,\b)}+kz)^{-1}$, we need to evaluate all other factors of
the product at $z=-\a^*U_{j_{+}(\a,\b)}/k=-\chi_{\a,\b}/k$. 
Recalling from 2.4 that 
\[ \a^*U_j - kU_j(d_{\a,\b})\frac{\chi_{\a,\b}}{k}=\b^*U_j,\] 
we find:
\footnote{Let us remind ourselves that we are using the analytic continuation 
convention $\prod_{m=1}^n:=\prod_{m=-\infty}^n /\prod_{m=-\infty}^0$.} 
\begin{align*} \text{For}&\ j\notin \a\cup \b, \ \ 
\prod_{m=1}^{U_j(d)+\a^*U_j(D)}(\a^*U_j-m\frac{\chi_{\a,\b}}{k}) = \\
& \prod_{m=1}^{kU_j(d_{\a,\b})} (\a^*U_j-m\frac{\chi_{\a,\b}}{k})\times \  
\prod_{m=1}^{U_j(d-kd_{\a,\b})+\a^*U_j(D)} 
(\b^*U_j-m\frac{\chi_{\a,\b}}{k}), \\
\text{for}&\ j=j_{-}(\a,\b),\ \ \prod_{m=1}^{U_j(d)} 
(-m\frac{\chi_{\a,\b}}{k}) = \\ & \prod_{m=1}^{U_j(kd_{\a,\b})}
(\a^*U_j-m\frac{\chi_{\a,\b}}{k})\ \times \ 
\prod_{m=1}^{U_j(d-kd_{\a,\b})+\a^*U_j(D)} 
(\b^*U_j-m\frac{\chi_{\a,\b}}{k})   \\
\text{for}&\ j=j_{+}(\a,\b),\ \ 
\prod_{m=1,\neq k}^{U_j(d)+\a^*U_j(D)}
(\a^*U_j-m\frac{\chi_{\a,\b}}{k}) = \\
& \prod_{m=1}^{k-1} (\a^*U_{j_{+}(\a,\b)}-m\frac{\chi_{\a,\b}}{k})\ \times 
\prod_{m=1}^{U_j(d-kd_{\a,\b})+\a^*U_j(D)-\b^*U_j(D)} 
(-m\frac{\chi_{\a,\b}}{k}), \\
 \text{for}&\ j\in \a \cap \b, \ \ \prod_{m=1}^{U_j(d)} 
(-m\frac{\chi_{\a,\b}}{k})= 1\ \times \ \prod_{m=1}^{U_j(d-kd_{\a,\b})
+\a^*U_j(D)-\b^*U_j(D)} 
(-m\frac{\chi_{\a,\b}}{k}), \\
& (Q^Dq^{P^{\a}(D)}) q^d = q^{kd_{\a,\b}} \times (Q^D q^{P^{\b}(D)})\ 
q^{d-kd_{\a,\b}+P^{\a}(D)-P^{\b}(D)}, \\
\exp & \left( -\frac{P^{\a}tk}{\chi_{\a,\b}}\right) 
\exp (dt) \exp (P^{\a}(D)t )=\\ & 1 \ \times \  
\exp \left(-\frac{P^{\b}tk}{\chi_{\a,\b}}\right)  
\exp \left[\left(d-kd_{\a,\b}+P^{\a}(D)-P^{\b}(D)\right)t\right] 
\exp (P^{\b}(D)t).
\end{align*}
In the last equality we use $(P^{\a}-P^{\b})/\chi_{\a,\b}=d_{\a,\b}$ from 2.4.

Factors on the R.H.S. which come before the multiplication 
sign ``$\times$'' form the recursion coefficients 
$q^{kd_{\a,\b}}Euler^{-1}_T({\mathcal N}_{\a,\b}(k))$. Factors which
come after the multiplication sign form the term of the series $\J^{\b}$
evaluated at $z=-\chi_{\a,\b}/k$ and with the summation index $d$ replaced 
with $d-kd_{\a,\b}+P^{\a}(D)-P^{\b}(D)$. Reversing this change in the 
summation index, we conclude that
\[ \Res_{z=-\frac{\chi_{\a,\b}}{k}}\ \J^{\a}(z)\ dkz = \frac{q^{kd_{\a,\b}}}
{Euler_T({\mathcal N}_{\a,\b}(k))} \ \J^{\b}(-\frac{\chi_{\a,\b}}{k}),\]
as required.

\section{Asymptotics}

{\bf 4.1. Stationary phase asymptotics.} We discuss here basic properties of
complex oscillating integrals 
\[ \int e^{f(x)/z} a(x) dx \] 
and their asymptotics as $z\to 0$. For simplicity of notation we assume all
integrals one-dimensional. Generalizations to higher dimensions are 
straightforward and are left to the reader. 

Let $x=0$ be a non-degenerate critical point of the {\em phase function} $f$,
i.e. $f(x)=f(0)-x^2/2\sigma^2 + \a x^3 +\b x^4+\dots $, and let $a(x)=
a_0+a_1x+a^2x^2+\dots$. The {\em stationary phase asymptotics} of the 
oscillating integral assumes the form:
\[ \int e^{f(x)/z} a(x) dx \sim \sqrt{2\pi z}\ \sigma\ e^{f(0)/z}
\sum_{k\geq 0} A_k z^k,\]  
where the coefficients $A_k$ are obtained by the following procedure. 
Make the change $x=\sqrt{z}\ y$, and replace a fixed integration 
interval $-m\leq x \leq m$ with the infinite interval $-\infty< y < \infty$
(to which $[-\frac{m}{\sqrt{z}} , \frac{m}{\sqrt{z}} ]$ tends as $z\to 0$):    
\[ \sqrt{z}\ e^{f(0)/z} \int_{-\infty}^{\infty} e^{-y^2/2\sigma^2} 
\exp (\a \sqrt{z}\ y^3+\b zy^4+\dots )\ (a_0+a_1\sqrt{z}\ y+a_2zy^2+\dots )
\ dy.\]
Expanding the integrand as a power series in $\sqrt{z}$ and evaluating
momenta of the Gaussian distribution
\[  \int_{-\infty}^{\infty} e^{-y^2/2\sigma^2} y^n dy = \left\{ 
\begin{array}{cl}
0 & \text{if $n$ is odd}, \\ \sqrt{2\pi}\ \sigma^{n+1} (n-1)!! & 
\text{if $n$ is even,} \end{array} \right. \]
we obtain the required asymptotical expansion.
Note that the sign of
$\sigma:=1/\sqrt{-f''(0)}$ depends on the choice of a branch of 
the square root, but the values of the {\em asymptotical coefficients} 
$A_k$ do not.
In this construction, the {\em amplitude}\ \ $a$ may depend formally
on $z$. Also, the phase function $f$ and/or the amplitude $a$ 
may depend on additional parameters, in which case the 
{\em critical value} $f(0)$, {\em Hessian} $-\sigma^{-2}$, and 
asymptotical coefficients do too. 
The following (rather obvious) 
proposition also allows for such parametric dependence. 
     
{\tt Proposition 1.} {\em Suppose that the 1-form in the integrand 
of an oscillating 
integral is the total Lie derivative along a vector field $v(x)\p/\p x$. 
Then the stationary phase asymptotics of this integral is trivial:
\[ \int d \left( e^{f(x)/z} a(x) v(x) \right) \sim 0 .\] }

{\em Proof.} After the change $x=\sqrt{z}\ y$ and series expansion, 
we arrive at the sequence of integrals
\[ \int_{-\infty}^{\infty} d \left( e^{-y^2/2\sigma^2} y^n \right) =  
\int_{-\infty}^{\infty} e^{-y^2/2\sigma^2} \left( ny^{n-1}-y^{n+1}\sigma^{-2}
\right) dy = 0.\]
The proposition follows from the (obvious) fact that the zero answer is 
obtained by substituting respective values of momenta of the Gaussian 
distribution (in lieu of actual integration).  Taking $v=1$ we obtain:

{\tt Corollary 1} (integration by parts). {\em The following oscillating 
integrals have the same stationary phase asymptotics:
\[ \int e^{f/z} f'g\ dx\ \ \ \text{and} \ \ \ -z \int e^{f/z} g'\ dx.\]} 
\indent The following two corollaries follow from their infinitesimal version 
established by Proposition 1.
 
{\tt Corollary 2.} {\em If a one-parameter family of oscillating integrals 
is obtained from each other by the flow of a vector field, then the  
stationary phase asymptotics does not depend on the parameter.}

{\tt Corollary 3.} {\em Stationary phase asymptotics of an oscillating integral
does not change under (formal or analytic) change of variables in the integral 
in a neighborhood of the non-degenerate critical point.}    
   
{\tt Proposition 2.} {\em The asymptotics of the derivative of an 
oscillating integral with respect to a parameter is obtained by 
differentiating the asymptotics of the integral.}

{\em Proof.} Let
\[ \int e^{f(x,\eps )/z} a(x, \eps) dx \sim \sqrt{2\pi z}\ \sigma (\eps)\  
e^{f(x_{cr}(\eps))/z} \sum_k A_k(\eps) z^k.\]
The RHS is obtained by evaluating momenta of Gaussian distributions after
the change 
\[ x=x_{cr}(\eps)+\sqrt{z}\ \sigma (\eps) y,\]
 where $x_{cr}(\eps)$ is the non-degenerate critical point of the phase 
function 
depending on the parameter, and $1/\sigma^2=f''(x_{cr}(\eps),\eps)$.  
Applying  $z\p/\p \eps$ to the RHS is equivalent to
differentiating the integrand termwise {\em after} the change of variables.  
On the other hand, we have:
\[ z\frac{\p}{\p \eps} \int e^{f/z} a dx = \int e^{f/z} \left(
\frac{\p f}{\p \eps} a + z\frac{\p a}{\p \eps} \right) dx.\]
Since the RHS has the same phase function as the original integral, 
the asymptotics of the derivative integral is obtained by applying 
the same operations: the change of variables, expansion of the integrand 
into a series, and evaluation of momenta, {\em preceded} however 
by the differentiation of the initial integrand. 
The change $(y,\eps) \mapsto (x (y,\eps),\eps)$ transforms
$\frac{\p}{\p \eps}$  into 
\[ \frac{\p}{\p \eps} + \frac{\p x}{\p \eps} 
\frac{\p}{\p x}.\] 
The difference with $\p/\p \eps$ is a vector field $(\p x/ \p \eps) \p/\p x$ 
(depending on
the parameter $\eps$). Thus the required independence of the asymptotics of 
the order of the operations follows from Proposition 1.        

{\tt Corollary.} {\em Given an oscillating integral 
\[ \int e^{f(x,\eps)/z} a(x,\eps) dx \sim \sqrt{2\pi z}\ \sigma (\eps)\  
e^{f(x_{cr}(\eps),\eps)/z} \sum_{k\geq 0} A_k (\eps) z^k,\]
depending formally on $\eps$, and
such that $f(x,0)$ has a non-degenerate critical point $x_{cr}(0)$,
consider it as an oscillating integral with the phase function $f(x,0)$ 
and the amplitude depending formally on $z$ and $\mu =\eps/z$: 
\[ \int e^{f(x,0)/z} \exp \left( \frac{f(x,z\mu)-f(x,0)}{z} \right) a(x,z\mu) 
dx.\]  
Then the asymptotics of the latter oscillating integral coincides
with the asymptotics of the former one at $\eps =z\mu$:  
\[ \sqrt{2\pi z}\ \sigma (z\mu) e^{f(x_{cr}(z\mu), z\mu)/z} \sum_{k\geq 0}
A_k (z\mu) z^k.\]}     

Indeed, it suffices to check that for each $n\geq 0$, the $n$th 
derivatives $(\p/\p \mu)^n$ of the two asymptotics coincide at $\mu=0$.
This follows by iterative application of Proposition 2 to the derivation
$z\p/\p \eps = \p/\p \mu$.  

{\bf 4.2. D-modules generated by J-functions.} Let $\LL \subset \H$ be an 
overruled Lagrangian cone in a symplectic loop space $(\H, \Omega)$, and let
$J$ be its J-function.  Tangent
spaces to $\LL$ vary in a family $H\ni \tau \to T_{\tau}$, where 
$J(-z,\tau)\in \LL$ is taken for the application point of $T_{\tau}$.
The J-function satisfies a system of 2nd order PDE:
\[  z\frac{\p^2 J(z,\tau)}{\p \tau^{\a}\p \tau^{\b}} = \sum_{\c} F_{\a\b}^{\c}
(\tau) \frac{\p J(z,\tau) }{\p \tau^{\c}},\]
where $\tau =\sum \tau^{\a} \phi_{\a}$ is a coordinate system on $H$. 

Indeed, for any family $\tau \to I(-z,\tau) \in \LL$ transverse to the ruling 
subspaces $zT_{\tau} \subset \LL$, the derivatives $\p I(-z,\tau)/\p \tau^{\b}$
form a basis of $T_{\tau}$ as a $\QQ [z]$-module, while 
$z \p I(-z,\tau)/\p \tau^{\b}$ lie in $zT_{\tau}\subset \LL$. Therefore the 2nd
derivatives $z\p^2 I(-z,\tau)/\p \tau^{\a}\p\tau^{\b}$ lie in $T_{\tau}$ and
are expressible as linear combinations of the basis, i.e.
\[  z\frac{\p^2 I(z,\tau)}{\p \tau^{\a}\p \tau^{\b}} = \sum_{\c} A_{\a\b}^{\c}
(z, \tau) \frac{\p I(z,\tau) }{\p \tau^{\c}}\]
where $A_{\a\b}^{\c}$ are suitable coefficients polynomial in $z$. When $I$
is the J-function, the LHS lies in $z\H_{-}$, while 
$\p I(-z,\tau)/\p \tau^{\c}$ form a basis of the quotient space 
$z\H_{-}/\H_{-}$. Thus the RHS lies in $z\H_{-}$ only if the coefficients
$A_{\a\b}^{\c}$ do not depend on $z$.\footnote{This line of reasoning is due 
to S. Barannikov \cite{Barannikov}.} 
 
In fact the coefficients $F_{\a\b}^{\c}(\tau)$ (and more generally, the values
$A_{\a\b}^{\c}(0,\tau)$) are structure constants of the
{\em quantum cup-product} $\bullet$:
\[ \phi_{\a}\bullet \phi_{\b} = \sum_{\c} F_{\a\b}^{\c}(\tau) \phi_{\c}.\]
Together with the pairing $(\cdot,\cdot)$, it provides the tangent spaces
$T_{\tau}H$ with the structure of a Frobenius algebra. In Gromov--Witten
theory,
\[ F_{\a\b}^{\c}=\sum_{n,d}\frac{Q^d}{n!}\int_{[M_{0,n+3,d}]}
\prod_{i=1}^n \ev_i^*(\tau) \ \ev_{n+1}^*(\phi_{\a}) \ev_{n+2}^*(\phi_{\b}) 
\ev_{n+3}^*(\phi^{\c}),\]
where $\{ \phi^{\c} \}$ is the basis of $H$ Poincare-dual to $\{ \phi_{\c} \}$.
  
Using the quantum cup-product, we can rewrite the PDE system for the J-function
in a more invariant form: 
\[ \forall v,w\in H,\ \ z\p_v \p_w J = \p_{v\bullet w} J.\] 
These equations can be considered as defining relations of the D-module
generated by the J-function. They allow one to represent 2nd derivatives 
of $J$ as linear combinations of first derivatives. Note that $v\bullet w$ 
depends on $\tau$. As a result, further differentiations of these equations 
contain terms involving derivatives of $v\bullet w$. However such terms come
with an extra $z$. Arguing inductively, we conclude:

{\tt Proposition 3} (\cite{Coates-Givental}). 
{\em For any $v_1, \dots, v_m \in H$, the higher directional derivatives 
$(z\p_{v_1})\cdots (z\p_{v_m}) J(z,\tau)$ can be expressed as linear
combinations of 1st derivatives $z\p_{\phi^{\c}} J(z,\tau)$ with coefficients
which are functions of $(z,\tau)$ polynomial in $z$. Modulo $(z)$, the 
direction vector of this linear combination coincides with the quantum 
cup-product $v_1\bullet \cdots \bullet v_m$, i.e.
\[ (z\p_{v_1})\cdots (z\p_{v_m}) J = z\p_{v_1\bullet \cdots \bullet v_m} J 
+ o(z).\]}

{\bf 4.3. Action of pseudo-differential symbols on J-functions.} 
We describe here in a general form a key argument from the proof
of Quantum Lefschetz Theorem found in \cite{Coates-Givental}. 

Let $(p_1, \dots, p_n)$ be coordinates on the space $H^*$ corresponding to the 
basis $( \phi_1, \dots,\phi_n)$ of $H$. Let $\Phi (p_1,\dots,p_n)$ be 
a polynomial. We consider it as the symbol of a differential operator 
$\Phi (z\p_{\phi_1}, \dots, z\p_{\phi_n})$ with constant coefficients.  

{\tt Lemma.} {\em Adjoin a formal parameter $\nu$ 
to the ground ring $\QQ$, and consider the overruled Lagrangian cone $\LL$ 
completed in the $\nu$-adic topology of $\QQ [[\nu]]$. Then 
\[ e^{-\nu \Phi (-z\p)/z} J (-z,\tau) \in \LL.\] }   
  
{\em Proof.} Let us assume first that $\Phi (0)=0$.  According to
Proposition 3, the action of the high order differential operator
$\Phi (z\p)$ on the J-function can be, in the quasi-classical
approximation, replaced with the action of the 1st order operator
$z\p_{\Phi (p\bullet)}$ where $\Phi (p\bullet) = \Phi (\phi_1\bullet,
\dots, \phi_n\bullet ) 1$ is the value of $\Phi$ computed in the quantum
cohomology algebra. Therefore $e^{\nu \Phi (z\p)/z} J$ can be written as
\[ \left[ 1+\nu \sum_{\c=1}^n a_{\c}(z,\tau; \nu) z\p_{\phi_\c} \right] \ 
e^{\nu \p_{\Phi (p\bullet)}} J(z,\tau),\]
where the coefficients $a_{\c}$ reduced modulo any power of $\nu$ 
are polynomial in $z$. By Taylor's formula, 
\[ I:=e^{\nu \p_{\Phi (p\bullet)}} J(z,\tau)=J(z,\tau+\nu\Phi (p\bullet)). \]
The point $I(-z,\tau)$ lies in $\LL$ since it is the value of the J-function, 
only at a shifted point. Thus the whole expression is obtained by adding to
this value a linear combination of the derivatives $z \p_{\phi_{\a}} I$ 
which lie in $zT_I(-z,\tau) \subset \LL$, so that the whole sum lies in $\LL$. 

In the case when $\Phi (0)\neq 0$, it suffices to add that 
\[ e^{\nu \Phi (0)/z} J(z,\tau)=J(z, \tau + \nu \Phi (0))\]
due to the {\em string equation} $z\p_1 J = J$.     

{\tt Corollary.} {\em The conclusion of the Lemma remains true even if
the differential symbol $\Phi (p, z, \nu)$ is allowed to depend on 
$z$ and $\nu$ (provided that modulo any power of $\nu$ it is polynomial in $z$).} 

{\em Remark.} In Proposition $3$, and hence in the results of $4.3$, 
one can replace the J-function by any function $\tau \mapsto I(-z,\tau) 
\in \LL$ transverse to the ruling spaces $zT_{\tau}$.     

\section{Mirrors}

{\bf 5.1. Mirrors of toric manifolds.} In equivariant Gromov--Witten theory
of toric manifolds, the {\em mirror} of the toric manifold $X$ (see $1.2$ 
for notations) is defined as the following oscillating integral
\[ \I (z,qe^t, \l):= \int 
 e^{\sum_{j=1}^N (x_j+\l_j \ln x_j)/z} 
\frac{d\ln x_1\w \cdots \w d\ln x_N}{d\ln q_1e^{t_1}\w \cdots \w 
d\ln q_Ke^{t_K}}\]
over suitable cycles in subvarieties of $\C^N$ given by the equations: 
\[  \prod_{j=1}^N x_j^{m_{ij}}=q_ie^{t_i}, \ \  i=1,...,K.\]
To a fixed point $\a =(j_1,\dots, j_K)\in X^T$, one associates a cycle
$C_{\a}$ of integration which is $\R_{+}^{N-K}$ in the chart 
$\{ x_j | j\notin \a\}$. 
On this cycle, the variables $x_j$ with $j\in \a$ can be expressed 
via the above relations in terms of the coordinates $x_j, j\notin \a$. 
Put
\[ \I_{\a}(z,t,q,\l):= \int_{C_{\a}}  
e^{\sum_{j=1}^N (x_j+\l_j \ln x_j)/z} {\bigwedge}_{j \notin \a} d\ln x_j.\]

{\tt Theorem 3.} {\em Let $J(z,\tau, Q)$ be the J-function of 
the base $B$ of the toric fibration $E\to B$ with the fiber $X$. Then  
\[ q^{-P^{\a}/z}\I_{\a}(z,t,q, z\p_{\L})J(z,\tau, Q)=\a^*I_E(z,t,\tau,q,Q) 
\prod_{j\notin \a}\!\int_{0}^{\infty}\!e^{(x-\a^*U_j\ln x)/z}d\ln x.  \]}

{\em Proof.} As it was mentioned in 1.1, the dependence of the genus-0 
descendant potential on Novikov's variables is 
governed by {\em divisor equations}. 
They can be stated in terms of the cone $\LL_B$ as follows
\cite{Coates-Givental, Givental_symplectic}. 
Novikov's variables $Q_i$ represent degrees (of holomorphic curves) which form 
a basis in $H_2(B,\Q)$. Let $\{ \rho_i \}$ denote the dual basis in 
$H^2(B,\Q)$. The linear operator $f\mapsto \rho_i f/z$ lies in the Lie algebra
of the twisted loop group and thus defines in the symplectic 
loop space $(\H, \Omega)$ a linear Hamiltonian vector field which we denote 
$\rho_i/z$. Then $\LL_B$, considered as a family of Lagrangian cones depending 
on $Q$, is invariant under the flows of the vector fields 
$Q_i\p/\p Q_i - \rho_i/z$. 

The divisor equations give rise to the following symmetries of the 
J-function. Let $J(z,\tau,Q)=\sum_D Q^D J_D(z,\tau)$, and 
$\rho \in H^2(B,\Z)$. Then 
\[ J(z,\tau+t\rho,Q)=e^{\rho t/z} \sum_D Q^D e^{\rho(D)t} J_D(z,\tau),\]
where $\rho (D)$ is the value of the cohomology class $\rho$ 
on the homology class $D$.

Thus, we have: 
\[ J(z,\tau+\sum_j \L_j \ln x_j, Q) 
=e^{\sum_j \L_j\ln x_j/z} \sum_D J_D(z,\tau)
Q^D \prod_j x_j^{\L_j(D)}.\] 
Therefore
\begin{align*} \I_{\a}(z,t, q, z\p_{\L}) & J(z,\tau, Q) =
\int_{C_{\a}} e^{\sum_{j=1}^N (x_j+z\p_{\L_j} \ln x_j)/z} J(z,\tau) \  
{\bigwedge}_{j \notin \a} d\ln x_j = \\
&\int_{C_{\a}} 
e^{\sum_{j=1}^N x_j/z} J(z, \tau+\sum_{j=1}^N \L_j \ln x_j) \ 
{\bigwedge}_{j \notin \a} d\ln x_j = \\
\sum_D J_D(z,\tau) & Q^D\ \int_{C_{\a}} 
e^{\sum_{j=1}^N (x_j+\L_j \ln x_j)/z } \prod_{j=1}^N 
x_j^{\L_j(D)} {\bigwedge}_{j \notin \a} d\ln x_j. \end{align*}  
On the other hand, relations between $x_j$ can be written in a 
more general form:
\[ \forall d\in \Z^K,  \ \ \ (qe^t)^d=\prod_j x_j^{U_j(d)},\]
since $U_j(d)=\sum_i d_im_{ij}$. The map $\Z^K\ni d \mapsto \{ U_j(d), \ j\in 
\a\} \in \Z^K$ is an isomorphism of lattices (this is a necessary condition 
for the toric variety $X$ to be non-singular at the fixed point $\a$). 
Therefore
\[ e^{\sum_{j\in \a} x_j/z} = \sum_{\{ d\ |\ U_j(d)\geq 0 \ \forall j\in \a\}}
\frac{(qe^t)^d}{\prod_{j\in \a} z^{U_j(d)} U_j(d)! \ 
\prod_{j\notin \a} x_j^{U_j(d)}}.\] 
Furthermore, from $\L_j=\sum_i P_i^{\a}m_{ij}-\a^*U_j$ (see $2.1$), we find:
\begin{align*} 
e^{\sum_{j=1}^N \L_j\ln x_j/z} & =\prod_{j=1}^N x_j^{-\a^*U_j/z}\prod_{i=1}^K
\left[ \prod_{j=1}^N x_j^{m_{ij}} \right]^{P_i^{\a}/z} = 
\prod_{j\notin \a} x_j^{-\a^*U_j/z} \prod_{i=1}^K (q_ie^{t_i})^{P_i^{\a}/z},\\
\prod_{j=1}^N x_j^{\L_j(D)}=&\prod_{j=1}^N x_j^{-\a^*U_j(D)} \prod_{i=1}^K
\left[ \prod_{j=1}^N x_j^{m_{ij}} \right]^{P_i^{\a}(D)} =\prod_{j\notin \a}
x_j^{-\a^*U_j(D)} \prod_{i=1}^K (q_ie^{t_i})^{P^{\a}_i(D)}. \end{align*}
Using this we rearrange the integrand to obtain:  
\begin{align*} q^{-P^{\a}/z} & \I_{\a}(z,t,q, z\p_{\L}) J(z,\tau, Q) = 
 e^{P^{\a}t/z}\sum_D J_D(z,\tau) Q^D (qe^t)^{P^{\a}(D)} \times \\
\sum_{d\in \Z^K} & \frac{(qe^t)^d}{\prod_{j\in\a} \prod_{m=1}^{U_j(d)}(mz)} 
 \prod_{j\notin \a} \int_{0}^{\infty} e^{x/z} 
x^{-\a^*U_j/z-U_j(d)-\a^*U_j(D)-1}\ dx.\end{align*}
Integrating by parts $U_j(d)+\a^*U_j(D)$ times, and making assumptions 
about the values of $z$ and $\a^*U_j$ which would guarantee that the 
integrand vanishes at $x=0$ and $x=\infty$, we find:
\[ \int_{0}^{\infty} e^{x/z} 
x^{-\a^*U_j/z-U_j(d)-\a^*U_j(D)-1}\ dx = \frac{\int_0^{\infty} 
e^{(x-\a^*U_j \ln x)/z} d\ln x}
{\prod_{m=1}^{U_j(d)+\a^*U_j(D)}(\a^*U_j+mz)}.\]
To complete the proof, substitute this into the previous formula, 
and compare the result with the expression for $\a^*I_E$ from Section 3:
 \begin{align*} \a^*I_E=e^{P^{\a}t/z} \sum_D\sum_{d \in \Z^K} 
\frac{J_D(z,\tau) (Q^Dq^{P^{\a}(D)}) q^d e^{dt} e^{P^{\a}(D)t}}
{\prod_{j\in \a} \prod_{m=1}^{U_j(d)} (mz)\ \prod_{j\notin \a} 
\prod_{m=1}^{U_j(d)+\a^*U_j(D)} (\a^*U_j+mz)}.
\end{align*}  

{\em Remark.} The assumptions about $\operatorname{Re} z$ and 
$\operatorname{Re} \L_j$, which guarantee
convergence of the integrals and vanishing of the finite terms that come out
of integration by parts, may differ for different terms of the series. 
The theorem should be understood therefore as the identity between 
coefficients of $(Q,q)$-series. In the next corollary about {\em asymptotics} 
of the integrals, convergence of the integrals is not required, and 
according to Corollary 1 of Proposition 1, integrations by parts does not
generate finite terms. 

{\tt Corollary.} {\em Let\ $\hat{\I}_{\a} (z,t,q, \l)$\ and\ 
$\hat{\Gamma} (z,\nu)$\  denote stationary phase asymptotics 
of the oscillating integrals $\I_{\a}$ and 
$\int_0^{\infty} e^{(-x+\nu\ln x)/z} d\ln x$ respectively. \linebreak 
Then
\[ q^{-P^{\a}/z} \hat{\I}_{\a}(z,t,q,z\p_{\L}) J(z,\tau, Q) = 
\a^*I_E(z,t,\tau,q,Q) \prod_{j\notin \a} \hat{\Gamma} (-z,\a^*U_j).\]}    

{\em Proof.} On the LHS, we have
\begin{align*} \hat{\I}_{\a}(z,t,q,z\p_{\L}) J(z,\tau, Q) &= \sum_D Q^D 
\hat{\I}_{\a}(z,t,q,z\p_{\L}) J_D(z,\tau) \\ &= \sum_D Q^D J_D(z,\tau) 
\hat{\I}_{\a}(z,t,q, \L+z\L(D)),\end{align*}
since $z\p_{\L_j} J_D = (\L_j+z\L_j(D)) J_D$ due to the divisor equation. 
The factor $\hat{\I}_{\a}(z,t,q, \L+z\L(D))$ is the stationary phase 
asymptotics of the integral 
\[ \int_{C_{\a}} e^{\sum_{j=1}^N (x_j+\L_j \ln x_j)/z } \prod_{j=1}^N 
x_j^{\L_j(D)} {\bigwedge}_{j \notin \a} d\ln x_j, \] 
which depends on the parameters $q$ (as well as $t$ and $\L$).
According to Corollary of Proposition 2, such asymptotics of a single 
oscillating integral depending on parameters can be replaced with
a suitable $q$-series of asymptotics of oscillating integrals. 
Following the steps in the proof of Theorem 3 and applying Corollary 1 of
Proposition 1 to justify integration by parts, 
we arrive at the expression on the RHS.

{\bf 5.2. The Quantum Riemann--Roch theorem.} We have:
\[ \hat{\Gamma}(z,\nu) = \sqrt{\frac{2\pi z}{\nu}} \ 
\exp\left\{\frac{\nu \ln \nu-\nu}{z} +\sum_{m=1}^{\infty} 
\frac{B_{2m}}{2m (2m-1)}\left(\frac{z}{\nu}\right)^{2m-1}\right\}.\]
Here $B_{2m}$ are Bernoulli numbers, and the equality follows from the
well-known asymptotics of the logarithm of the Gamma-function 
$\Gamma (\nu/z)$.

Let $\LL$ and $\LL^{tw}$ be the overruled Lagrangian cones respectively:
of genus 0 Gromov--Witten theory of a target manifold $M$, and of such 
a theory {\em twisted} (in the sense of $2.3$) by a line bundle over $M$ 
with the equivariant 1st Chern class $\nu$. The cone $\LL$ lies in the
symplectic loop space $(\H, \Omega)$ based on the Poincare pairing $
(a,b)=\int_M ab$, while $\LL^{tw}$ lies in $(\H, \Omega^{tw})$ based on
$(a,b)^{tw}=\int_M ab/\nu$. The linear map $(\H, \Omega^{tw}) \to 
(\H, \Omega)$ defined by $f \mapsto f/\sqrt{\nu}$ is a symplectomorphism. 

{\tt Theorem} (\cite{Coates-Givental}).  
\[ \LL = \frac{\hat{\Gamma}(z,\nu)}{\sqrt{2\pi z}} \LL^{tw}.\]

{\bf 5.3. Completing the proof of Theorem 1.} We need to show that 
$\a^*I_{E} (-z)$ lies in the overruled Lagrangian cone $\LL^{\a}$ corresponding
to the genus-0 Gromov--Witten theory of $B$ twisted by the normal bundle
$N^{\a}$ of $\a(B)\subset E$. Due to the Quantum Riemann--Roch 
Theorem, it suffices to prove that, equivalently, 
\[ \prod_{j\notin \a} \frac{\hat{\Gamma}(z,\a^*U_j)}{\sqrt{2\pi z}} 
\ \a^*I_E (-z) \]
lies in the overruled Lagrangian cone of the ``untwisted'' theory.
Note that in this formula, Novikov's variables still occur in the form
$Q^Dq^{P^{\a}(D)}$ to account correctly for degrees of curves in $\a(B)\subset
E$ considered as curves in $E$. According to Corollary of Theorem 3, 
the above expression coincides with
\[ (2\pi z)^{-\dim N^{\a}/2}
q^{P^{\a}/z} \hat{\I}_{\a}(-z,t,q,-z\p_{\L}) \sum_D Q^D J_D(-z,\tau).\] 
We make the change $Q^D\mapsto Q^Dq^{-P^{\a}(D)}$ to restore the absolute 
meaning of Novikov's variables, and apply the divisor equation:
\[ q^{P^{\a}/z} \sum_D Q^D q^{-P^{\a}(D)} J_D(-z,\tau) = 
J(-z,\tau-P^{\a}\ln q, Q).\]
Then it remains to show that 
\[ (2\pi z)^{-\dim N^{\a}/2} 
\hat{\I}_{\a}(-z,t,q,-z\p_{\L}) J(-z,\tau-P^{\a}\ln q, Q) \ \in \ \LL_B. \]
Since $J(z,\tau, Q)$ is the J-function of $\LL_B$, 
this follows from the results of 4.3 about actions of pseudo-differential 
symbols on J-functions.


\newpage
               
\begin{thebibliography}{100}

\bibitem{Barannikov} S. Barannikov. {\em Quantum Periods- I.  Semi- infinite 
variations of Hodge structures.}  Internat. Math. Res. Notices  2001,  no. 23, 1243--1264.     
\bibitem{Brown et al} J. Brown, T. Coates, A. Givental, H.-H. Tseng. 
{\em Virasoro constraints for toric fibrations.} in progress.
\bibitem{Coates-Givental} T. Coates, A. Givental. {\em Quantum Riemann-Roch, Lefschetz and Serre.}  Ann. of Math. (2)  165  (2007),  no. 1, 15--53.  
\bibitem{Elezi} A. Elezi. {\em A mirror conjecture for projective bundles.}
Intern. Math. Res. Notices, 2005, No. 55., 3445--3458.  
\bibitem{Givental_MSRI} A. Givental. {\em Gromov--Witten invariants of
symplectic quotients.} A lecture at MSRI, March 2006.  
\bibitem{Givental_toric} A. Givental. {\em A mirror theorem for toric 
complete intersections.} Topological field theory, primitive forms and related topics (Kyoto, 1996),  141--175, Progr. Math., 160, Birkhäuser Boston, Boston, MA, 1998. 
\bibitem{Givental_symplectic} A. Givental. {\em Symplectic geometry of 
Frobenius structures.} Frobenius manifolds,  91--112, Aspects Math., E36, Vieweg, Wiesbaden, 2004.
\bibitem{Graber-Pandharipande} T. Graber, R. Pandharipande. {\em Localization
of virtual classes.} Invent. Math.  135  (1999),  no. 2, 487--518.   
\bibitem{Iritani} H. Iritani. {\em Quantum D- modules and generalized mirror 
transformations.}  Topology  47  (2008),  no. 4, 225--276.    
\bibitem{Kontsevich} M. Kontsevich. {\em Enumeration of rational curves
via tori actions.} The moduli space of curves (Texel Island, 1994),  335--368, Progr. Math., 129, Birkhäuser Boston, Boston, MA, 1995.


\end{thebibliography}
\enddocument